\nonstopmode
\documentclass{amsart}
\usepackage{amsmath, amssymb, amsthm, amscd}
\usepackage[all]{xy}
\newdir{ >}{{}*!/-10pt/@{>}}
\newdir^{ (}{{}*!/-5pt/@^{(}}
\newdir_{ (}{{}*!/-5pt/@_{(}}

\def\undersetbrace#1\to#2{\underbrace{#2}_{#1}}                                                          
\def\oversetbrace#1\to#2{\overbrace{#2}^{#1}}
\def\AMSunderset#1\to#2{\underset{#1}{#2}}
\def\AMSoverset#1\to#2{\overset{#1}{#2}}

\def\East#1#2{\overset{#1}{\longrightarrow}}

\swapnumbers

\newtheorem*{prop*}{Proposition}
\newtheorem{thm}[subsection]{Theorem}
\newtheorem*{thm*}{Theorem}

\newtheorem*{lem*}{Lemma}
\newtheorem{cor}[subsection]{Corollary}
\newtheorem*{cor*}{Corollary}
\newtheorem*{remark*}{Remark}

\newenvironment{demo}[1]{{\textit{#1.}}}{\par\smallskip}
\numberwithin{equation}{subsection}

\def\ign#1{}             
\def\o{\,\circ\,}
\def\X{\mathfrak X}
\def\al{\alpha}
\def\be{\beta}

\def\ep{\varepsilon}

\def\et{\eta}

\def\la{\lambda}

\def\si{\sigma}

\def\Ga{\Gamma}

\def\i{^{-1}}
\def\x{\times}
\def\p{\partial}
\let\on=\operatorname

\def\AMSonly#1{}

\def\Id{\on{Id}}

\def\Diff{{\on{Diff}}}

\title[A zoo of diffeomorphism groups] 
{A zoo of diffeomorphism groups on $\mathbb R^n$} 

\author{Peter W. Michor and David Mumford}

\address{
Peter W. Michor:
Fakult\"at f\"ur Mathematik, Universit\"at Wien,
Nordbergstrasse 15, A-1090 Wien, Austria}
\email{Peter.Michor@univie.ac.at}
\address{
David Mumford:
Division of Applied Mathematics, Brown University,
Box F, Providence, RI 02912, USA}
\email{David\_{}Mumford@brown.edu}
\date{\today}

\begin{document}
\begin{abstract}
We consider the groups $\Diff_{\mathcal B}(\mathbb R^n)$, $\Diff_{H^\infty}(\mathbb R^n)$, 
and $\Diff_{\mathcal S}(\mathbb R^n)$ of smooth diffeomorphisms on $\mathbb R^n$ which differ from 
the identity by a function which is in either $\mathcal B$ (bounded in all derivatives), 
$H^\infty = \bigcap_{k\ge 0}H^k$, or $\mathcal S$ (rapidly decreasing). We show that all these 
groups are smooth regular Lie groups.
\end{abstract}

\maketitle 

\section{Introduction}

The purpose of this article is to prove that the following groups of diffeomorphisms on 
$\mathbb R^n$ are regular (see \ref{reglie}) Lie groups:
\begin{itemize}
\item $\Diff_{\mathcal B}(\mathbb R^n)$, the group of all diffeomorphisms which differ from the 
       identity by a function which is bounded together with all derivatives separately; see \ref{diffb}.
\item $\Diff_{H^\infty}(\mathbb R^n)$, the group of all diffeomorphisms which differ from the 
       identity by a function in the intersection $H^\infty$ of all Sobolev spaces $H^k$ for 
       $k\in \mathbb N_{\ge 0}$; see \ref{diffh}.
\item $\Diff_{\mathcal S}(\mathbb R^n)$, the group of all diffeomorphisms which fall rapidly to the 
       identity; see \ref{diffs}.
\end{itemize} 
Since we are giving a kind of uniform proof, we also mention in \ref{diffc} the group $\Diff_{c}(\mathbb R^n)$ 
of all diffeomorphisms which differ from the 
identity only on a compact subset, where this result is known for many years, by \cite{Michor80II} 
and \cite{Michor80}. 
The groups $\Diff_{\mathcal S}(\mathbb R^n)$ and partly $\Diff_{H^\infty}(\mathbb R^n)$  have been used 
essentially in the papers 
\cite{Michor130},
\cite{Michor127},
\cite{Bauer2011c_preprint},
\cite{Bauer2012c_preprint},
\cite{Bauer2012d_preprint},
and \cite{Bauer2012c}.
In particular, $\Diff_{H^\infty}(\mathbb R^n)$ is essential if one wants to prove that the geodesic 
equation of a right Riemannian invariant metric is well-posed with the use of Sobolov space 
techniques. The regular Lie groups $\Diff_{\mathcal B}(E)$ and $\Diff_{\mathcal S}(E)$ have been 
treated, using single derivatives iteratively, in \cite{Walter12}, for a Banach space $E$. 
See \cite{Goldin04} for the role of diffeomorphism groups in quantum physics. 
Andreas Kriegl, Leonard Frerick, and Jochen Wengenroth helped with discussions and hints.

\section{Some words on smooth convenient calculus}

Traditional differential calculus works 
well for finite dimensional vector spaces and for Banach spaces. For 
more general locally convex spaces 
we sketch here the convenient approach as explained in 
\cite{FK88} and \cite{KM97}.
The main difficulty is that composition of 
linear mappings stops to be jointly continuous at the level of Banach 
spaces, for any compatible topology. 
We use the notation of \cite{KM97} and this is the
main reference for this section. 

\subsection{The $c^\infty$-topology}
Let $E$ be a 
locally convex vector space. A curve $c:\mathbb R\to E$ is called 
{\it smooth} or $C^\infty$ if all derivatives exist and are 
continuous - this is a concept without problems. Let 
$C^\infty(\mathbb R,E)$ be the space of smooth functions. It can be 
shown that the set $C^\infty(\mathbb R,E)$ does not depend on the locally convex 
topology of $E$, only on its associated bornology (system of bounded 
sets).
The final topologies with respect to the following sets of mappings into E coincide:
\begin{enumerate}
\item $C^\infty(\mathbb R,E)$.
\item The set of all Lipschitz curves 
(so that $\{\frac{c(t)-c(s)}{t-s}:t\neq s, |t|, |s|\le C\}$ 
is bounded in $E$, for each $C>0$). 
\item The set of injections $E_B\to E$ where $B$ runs through all bounded 
absolutely convex subsets in $E$, and where 
$E_B$ is the linear span of $B$ equipped with the Minkowski 
functional $\|x\|_B:= \inf\{\la>0:x\in\la B\}$.
\item The set of all Mackey-convergent sequences $x_n\to x$ 
(there exists a sequence 
$0<\la_n\nearrow\infty$ with $\la_n(x_n-x)$ bounded).
\end{enumerate}
This topology is called the $c^\infty$-topology on $E$ and we write 
$c^\infty E$ for the resulting topological space. In general 
(on the space $\mathcal{D}$ of test functions for example) it is finer 
than the given locally convex topology, it is not a vector space 
topology, since scalar multiplication is no longer jointly 
continuous. The finest among all locally convex topologies on $E$ 
which are coarser than $c^\infty E$ is the bornologification of the 
given locally convex topology. If $E$ is a Fr\'echet space, then 
$c^\infty E = E$. 

\subsection{Convenient vector spaces} 
A locally convex vector space 
$E$ is said to be a {\it convenient 
vector space} if one of the following equivalent
conditions is satisfied (called $c^\infty$-completeness):
\begin{enumerate}
\item For any $c\in C^\infty(\mathbb R,E)$ the (Riemann-) integral 
$\int_0^1c(t)dt$ exists in $E$.
\item Any Lipschitz curve in $E$ is locally Riemann integrable.
\item A curve $c:\mathbb R\to E$ is smooth if and only if $\la\o c$ is 
smooth for all $\la\in E^*$, where $E^*$ is the dual consisting 
of all continuous linear functionals on $E$. Equivalently, 
we may use the dual $E'$ consisting of 
all bounded linear functionals.
\item Any Mackey-Cauchy-sequence (i.\ e.\  $t_{nm}(x_n-x_m)\to 0$  
for some $t_{nm}\to \infty$ in $\mathbb R$) converges in $E$. 
This is visibly a mild completeness requirement.
\item If $B$ is bounded closed absolutely convex, then $E_B$ 
is a Banach space.
\item If $f:\mathbb R\to E$ is scalarwise $\on{Lip}^k$, then $f$ is 
$\on{Lip}^k$, for $k\ge0$.
\item If $f:\mathbb R\to E$ is scalarwise $C^\infty$ then $f$ is 
differentiable at 0.
\item If $f:\mathbb R\to E$ is scalarwise $C^\infty$ then $f$ is 
$C^\infty$.
\end{enumerate}
Here a mapping $f:\mathbb R\to E$ is called $\on{Lip}^k$ if all 
derivatives up to order $k$ exist and are Lipschitz, locally on 
$\mathbb R$. That $f$ is scalarwise $C^\infty$ means $\la\o f$ is $C^\infty$  
for all continuous linear functionals on $E$.

\subsection{Smooth mappings} 
Let $E$, $F$, and $G$ be convenient vector spaces, 
and let $U\subset E$ be $c^\infty$-open. 
A mapping $f:U\to F$ is called {\it smooth} or 
$C^\infty$, if $f\o c\in C^\infty(\mathbb R,F)$ for all 
$c\in C^\infty(\mathbb R,U)$.
{\it
The main properties of smooth calculus are the following.
\begin{enumerate}
\item For mappings on Fr\'echet spaces this notion of smoothness 
coincides with all other reasonable definitions. Even on 
$\mathbb R^2$ this is non-trivial.
\item Multilinear mappings are smooth if and only if they are 
bounded.
\item If $f:E\supseteq U\to F$ is smooth then the derivative 
$df:U\x E\to F$ is  
smooth, and also $df:U\to L(E,F)$ is smooth where $L(E,F)$ 
denotes the space of all bounded linear mappings with the 
topology of uniform convergence on bounded subsets.
\item The chain rule holds.
\item The space $C^\infty(U,F)$ is again a convenient vector space 
where the structure is given by the obvious injection
$$
C^\infty(U,F) \East{C^\infty(c,\ell)}{} 
\negthickspace\negthickspace\negthickspace\negthickspace\negthickspace
\negthickspace\negthickspace\negthickspace
\prod_{c\in C^\infty(\mathbb R,U), \ell\in F^*} 
\negthickspace\negthickspace\negthickspace\negthickspace\negthickspace
\negthickspace\negthickspace\negthickspace
C^\infty(\mathbb R,\mathbb R),
\quad f\mapsto (\ell\o f\o c)_{c,\ell},
$$
where $C^\infty(\mathbb R,\mathbb R)$ carries the topology of compact 
convergence in each derivative separately.
\item The exponential law holds: For $c^\infty$-open $V\subset F$, 
$$
C^\infty(U,C^\infty(V,G)) \cong C^\infty(U\x V, G)
$$
is a linear diffeomorphism of convenient vector spaces. 
\item A linear mapping $f:E\to C^\infty(V,G)$ is smooth (by \thetag{2} equivalent to bounded) if 
and only if $E \East{f}{} C^\infty(V,G) \East{\on{ev}_v}{} G$ is smooth 
for each $v\in V$. This is called the smooth uniform 
boundedness theorem \cite[5.26]{KM97}.
\item The following canonical mappings are smooth.
\begin{align*}
&\operatorname{ev}: C^\infty(E,F)\x E\to F,\quad 
\operatorname{ev}(f,x) = f(x)\\
&\operatorname{ins}: E\to C^\infty(F,E\x F),\quad
\operatorname{ins}(x)(y) = (x,y)\\
&(\quad)^\wedge :C^\infty(E,C^\infty(F,G))\to C^\infty(E\x F,G)\\
&(\quad)^\vee :C^\infty(E\x F,G)\to C^\infty(E,C^\infty(F,G))\\
&\operatorname{comp}:C^\infty(F,G)\x C^\infty(E,F)\to C^\infty(E,G)\\
&C^\infty(\quad,\quad):C^\infty(F,F_1)\x C^\infty(E_1,E)\to 
C^\infty(C^\infty(E,F),C^\infty(E_1,F_1))\\
&\qquad (f,g)\mapsto(h\mapsto f\o h\o g)\\
&\prod:\prod C^\infty(E_i,F_i)\to C^\infty(\prod E_i,\prod F_i)
\end{align*}
\end{enumerate}
}
Note that the conclusion of 
\thetag{6} is the starting point of the classical calculus of 
variations, where a smooth curve in a space of functions was assumed 
to be just a smooth function in one variable more. It is also the source of the name convenient 
calculus.
This and some other obvious properties already determine convenient calculus.
There are, however, smooth mappings which are not continuous. This is 
unavoidable and not so horrible as it might appear at first sight. 
For example the evaluation $E\x E^*\to\mathbb R$ is jointly continuous if 
and only if $E$ is normable, but it is always smooth. Clearly smooth 
mappings are continuous for the $c^\infty$-topology.
This ends our review of the standard results of convenient calculus.
But we will need more.

\begin{thm}\label{FrolicherKriegl}{\rm \cite[~4.1.19]{FK88}}
Let $c:\mathbb R\to E$ be a curve in a convenient vector space $E$. Let 
$\mathcal{V}\subset E'$ be a subset of bounded linear functionals such that 
the bornology of $E$ has a basis of $\sigma(E,\mathcal{V})$-closed sets. 
Then the following are equivalent:
\begin{enumerate}
\item $c$ is smooth
\item There exist locally bounded curves $c^{k}:\mathbb R\to E$ such that
      $\ell\o c$ is smooth $\mathbb R\to \mathbb R$ with $(\ell\o c)^{(k)}=\ell\o
      c^{k}$, for each $\ell\in\mathcal V$. 
\end{enumerate}
If $E$ is reflexive, then for any point separating subset
$\mathcal{V}\subset E'$ the bornology of $E$ has a basis of 
$\si(E,\mathcal{V})$-closed subsets, by {\rm \cite[~4.1.23]{FK88}}.
\end{thm}
This theorem is surprisingly strong: Note that $\mathcal V$ does not need to recognize bounded 
sets.

\subsection{Faa di Bruno formula}\label{Faa} 
Let $g\in C^\infty(\mathbb R^n,\mathbb R^k)$ and let $f\in C^\infty(\mathbb R^k,\mathbb R^l)$. 
Then the $p$-th deivative of $f\o g$ looks as follows where $\on{sym}_p$ denotes symmetrization of a $p$-linear 
mapping:
\begin{align*}
\frac{d^p(f\o g)(x)}{p!} &= \on{sym}_p\Big(
\sum_{j= 1}^p \sum_{\substack{\al\in \mathbb N_{>0}^j\\ \al_1+\dots+\al_j =p}}
\frac{d^jf(g(x))}{j!}\Big( 
\frac{d^{\al_1}g(x)}{\al_1!},\dots,
\frac{d^{\al_j}g(x)}{\al_j!}\Big)\Big)
\end{align*}
The one dimensional version is due to Fa\`a di Bruno \cite{FaaDiBruno1855}, the only
beatified mathematician. The formula is seen by composing the Taylor series.

\section{Groups of smooth diffeomorphisms}

\subsection{Model spaces for Lie groups of diffeomorphism}\label{modelspaces}
If we consider the group of all orientation preserving 
diffeomorphisms $\Diff(\mathbb R^n)$ of $\mathbb R^n$, it is not an open subset of 
$C^\infty(\mathbb R^n,\mathbb R^n)$ with the compact $C^\infty$-topology. 
So it is not a smooth manifold in the usual sense, but we may consider it as a Lie group in the 
cartesian closed category of  Fr\"olicher spaces, see \cite[section 23]{KM97} with the structure 
induced by the injection 
$f\mapsto (f,f\i)\in C^\infty(\mathbb R^n,\mathbb R^n)\x C^\infty(\mathbb R^n,\mathbb R^n)$.
Or one can use the theory of smooth manifolds based on smooth curves instead of charts from 
\cite{Michor84a}, \cite{Michor84b}, which agrees with the usual theory up to Banach manifolds.

We shall now describe regular Lie groups in $\Diff(\mathbb R^n)$ 
which are given by 
diffeomorphisms of the form $f = \Id + g$ 
where $g$ is in some specific convenient vector spaces of bounded functions in 
$C^\infty(\mathbb R^n,\mathbb R^n)$.
Now we discuss these spaces on $\mathbb R^n$, we describe the smooth curves in them, and we describe the 
corresponding groups.

\subsection{Regular Lie groups}\label{reglie}
We consider a  smooth Lie group $G$ with Lie algebra $\mathfrak g=T_eG$ modelled on convenient 
vector spaces. 
The notion of a regular Lie group is originally due to Omori and collaborators
(see \cite{OmoriMaedaYoshioka82}, \cite{OmoriMaedaYoshioka83}) for Fr\'echet Lie groups, was 
weakened and made more transparent by Milnor \cite{Milnor84} and carried over to convenient Lie 
groups in \cite{KM97r}, see also \cite[38.4]{KM97}.
A Lie group $G$
is called {\em regular}  if the following holds:
\begin{itemize}
\item 
For each smooth curve 
$X\in C^{\infty}(\mathbb R,\mathfrak g)$ there exists a curve 
$g\in C^{\infty}(\mathbb R,G)$ whose right logarithmic derivative is $X$, i.e.,
$$
\begin{cases} g(0) &= e \\
\p_t g(t) &= T_e(\mu^{g(t)})X(t) = X(t).g(t)
\end{cases} 
$$
where $\mu:G\x G\to G$ is  multiplication with $\mu(g,h)=g.h = \mu_g(h)= \mu^h(g)$.
The curve $g$ is uniquely determined by its initial value $g(0)$, if it
exists.
\item
Put $\on{evol}^r_G(X)=g(1)$ where $g$ is the unique solution required above. 
Then $\on{evol}^r_G: C^{\infty}(\mathbb R,\mathfrak g)\to G$ is required to be
$C^{\infty}$ also. 
\end{itemize}

\subsection{The group $\Diff_{\mathcal B}(\mathbb R^n)$}\label{diffb}
The space $\mathcal B(\mathbb R^n)$ (called $\mathcal D_{L^\infty}(\mathbb R^n)$ by L.~Schwartz 
\cite{Schwartz66}) consists of all smooth functions with all derivatives (separately) bounded.
It is a Fr\'echet space. 
By \cite{Vogt83}, the space $\mathcal B(\mathbb R^n)$ is linearly 
isomorphic to $\ell^\infty\hat\otimes\, \mathfrak s$ for any completed tensor product between the 
projective one and the injective one, where $\mathfrak s$ is the nuclear Fr\'echet space of rapidly 
decreasing real sequences. Thus $\mathcal B(\mathbb R^n)$ is not reflexive and not nuclear.
\newline
{\it The space $C^\infty(\mathbb R,\mathcal{B}(\mathbb R^n))$ of smooth
curves in $\mathcal{B}(\mathbb R^n)$ consists of all functions 
$c\in C^\infty(\mathbb R^{n+1},\mathbb R)$ satisfying the following 
property:
\begin{enumerate}
\item[$\bullet$]
For all $k\in \mathbb N_{\ge0}$, $\al\in \mathbb N_{\ge0}^n$ and each $t\in \mathbb R$ the expression
$\p_t^{k}\p^\al_x c(t,x)$  is uniformly bounded in $x\in \mathbb R^n$, locally
in $t$. 
\end{enumerate} 
}
To see this we use \ref{FrolicherKriegl} for the set $\{\on{ev}_x: x\in \mathbb R\}$ of point
evaluations in $\mathcal{B}(\mathbb R^n)$. 
Here $\p^\al_x = \frac{\p^{|\al|}}{\p x^\al}$ and $c^k(t)=\p_t^kf(t,\quad)$.
\newline
{\it $\Diff^+_{\mathcal B}(\mathbb R^n)=\bigl\{f=\Id+g: g\in\mathcal B(\mathbb R^n)^n, 
    \det(\mathbb I_n + dg)\ge \ep > 0 \bigr\}$ denotes the corresponding group}, 
		see theorem \ref{regularLie} below.

\subsection{The group $\Diff_{H^\infty}(\mathbb R^n)$}\label{diffh} 
The space 
$H^\infty(\mathbb R^n)=\bigcap_{k\ge 1}H^k(\mathbb R^n)$ is the intersection of all Sobolev 
spaces which is a reflexive Fr\'echet space.
It is called $\mathcal D_{L^2}(\mathbb R^n)$ by L.~Schwartz in \cite{Schwartz66}.
By \cite{Vogt83}, the space $H^{\infty}(\mathbb R^n)$ is linearly isomorphic to 
$\ell^2\hat\otimes\, \mathfrak s$. Thus it is not nuclear, not Schwartz, not Montel, but still smoothly 
paracompact.
\newline
{\it The space $C^\infty(\mathbb R,H^\infty(\mathbb R^n))$ of smooth
curves in $H^\infty(\mathbb R^n)$ consists of all functions 
$c\in C^\infty(\mathbb R^{n+1},\mathbb R)$ satisfying the following 
property:
\begin{enumerate}
\item[$\bullet$]
For all $k\in \mathbb N_{\ge0}$, $\al\in \mathbb N_{\ge0}^n$  the expression 
$\|\p_t^{k}\p^\al_xf(t,\quad)\|_{L^2(\mathbb R^n)}$ is locally bounded  
near each $t\in \mathbb R$. 
\end{enumerate} 
}
The proof is literally the same as for $\mathcal B(\mathbb R^n)$, noting that the point evaluations are 
continuous on each Sobolev space $H^k$ with $k>\frac n2$.
\newline
{\it $\Diff^+_{H^{\infty}}(\mathbb R^n)=\bigl\{f=\Id+g: g\in H^\infty(\mathbb R^n), 
		\det(\mathbb I_n + dg)>0\bigr\}$ 
		denotes the correponding group},
		see theorem \ref{regularLie} below.

\subsection{The group $\Diff_{\mathcal S}(\mathbb R^n)$}\label{diffs} 
The algebra $\mathcal S(\mathbb R^n)$  of rapidly decreasing functions is a reflexive nuclear Fr\'echet space.
\newline
{\it The space $C^\infty(\mathbb R,\mathcal S(\mathbb R^n))$ of smooth
curves in $\mathcal S(\mathbb R^n)$ consists of all functions 
$c\in C^\infty(\mathbb R^{n+1},\mathbb R)$ satisfying the following 
property:
\begin{enumerate}
\item[$\bullet$]
For all $k,m\in \mathbb N_{\ge0}$ and $\al\in \mathbb N_{\ge0}^n$,  the expression
$(1+|x|^2)^m\p_t^{k}\p^\al_xc(t,x)$ is uniformly bounded in $x\in \mathbb R^n$, locally uniformly bounded  
in $t\in \mathbb R$.
\end{enumerate} 
}
{\it $\Diff^+_{\mathcal S}(\mathbb R^n)=\bigl\{f=\Id+g: g\in \mathcal S(\mathbb R^n)^n, 
		\det(\mathbb I_n + dg)>0\bigr\}$ is the correponding group.} 
		
\subsection{The group  $\Diff_{c}(\mathbb R^n)$}\label{diffc} 
The algebra $C^\infty_c(\mathbb R^n)$ of all smooth functions with compact support is a nuclear 
(LF)-space. 
{\it The space $C^\infty(\mathbb R,C^\infty_c(\mathbb R^n))$ of smooth
curves in $C^\infty_c(\mathbb R^n)$ consists of all functions 
$f\in C^\infty(\mathbb R^{n+1},\mathbb R)$ satisfying the following 
property:
\begin{enumerate}
\item[$\bullet$]
For 
each compact interval $[a,b]$ in $\mathbb R$ there exists a compact subset $K\subset \mathbb R^n$
such that $f(t,x)=0$ for  $(t,x)\in [a,b]\x (\mathbb R^n\setminus K)$.
\end{enumerate} 
}
{\it $\Diff_c(\mathbb R^n)=\bigl\{f=\Id+g: g\in C^\infty_c(\mathbb R^n)^n, 
		\det(\mathbb I_n + dg)>0\bigr\}$ is the correponding group.}

\subsection{Ideal properties of function spaces}\label{ideal}			 
The function spaces are boundedly mapped into each other as follows:
$$\xymatrix{
C^\infty_c(\mathbb R^n) \ar[r]  & \mathcal S(\mathbb R^n)  \ar[r] & H^\infty(\mathbb R^n) \ar[r] 
 &  \mathcal B(\mathbb R^n)
}$$
and each space is a bounded locally convex algebra and a bounded $\mathcal B(\mathbb R^n)$-module.
Thus each space is an ideal in each larger space.

\begin{thm}\label{regularLie} 
The sets of diffeomorphisms $\Diff_c(\mathbb R^n)$, 
$\Diff_{\mathcal S}(\mathbb R^n)$,
$\Diff_{H^\infty}(\mathbb R^n)$, and 
$\Diff_{\mathcal B}(\mathbb R^n)$
are all smooth regular Lie groups in the sense of \ref{reglie}. 
We have the following smooth injective group homomorphisms
$$\xymatrix{
\Diff_c(\mathbb R^n) \ar[r] & \Diff_{\mathcal S}(\mathbb R^n) \ar[r] & \Diff_{H^\infty}(\mathbb R^n) \ar[r]  
&  \Diff_{\mathcal B}(\mathbb R^n)
}.$$
Each group is a normal subgroup in any other in which it is contained, in particular in 
$\on{Diff}_{\mathcal B}(\mathbb R^n)$.
\end{thm}

The case $\Diff_c(\mathbb R^n)$ is well-known, see for example 
\cite[43.1]{KM97}.
The 1-dimensional version 
$\Diff_{\mathcal S}(\mathbb R)$ was treated in \cite[6.4]{Michor109}.

\begin{proof} Let $\mathcal A$ denote any of $\mathcal B$, $H^\infty$, 
$\mathcal S$, or $c$, and let $\mathcal A(\mathbb R^n)$ denote the corresponding function space as 
described in \ref{diffb} - \ref{diffc}.
Let $f(x)= x+g(x)$ for $g\in \mathcal A(\mathbb R^n)^n$ 
with $\det(\mathbb I_n + dg)>0$ and for $x\in \mathbb R^n$.
We have to check that each $f$ as described is a diffeomorphism. By the inverse function theorem, 
$f$ is a locally a diffeomorphism everywhere. Thus the image of $f$ is open in $\mathbb R^n$.
We claim that it is also closed. So let $x_i\in \mathbb R^n$ with $f(x_i)=x_i + g(x_i)\to y_0$ in 
$\mathbb R^n$. Then $f(x_i)$ is a bounded sequence. 
Since $g\in \mathcal A(\mathbb R^n)\subset \mathcal B(\mathbb R^n)$, the $x_i$ 
also form a bounded sequence, thus contain a convergent subsequence. Without loss let $x_i\to x_0$ 
in $\mathbb R^n$. Then $f(x_i)\to f(x_0)=y_0$. Thus $f$ is surjective. 
This also shows that $f$ is a proper mapping (i.e., compact sets have compact inverse images under 
$f$). By \cite[17.2]{MichorH}, a proper surjective submersion is the projection of a smooth fiber 
bundle. In our case here $f$ has discrete fibers, so $f$ is a covering mapping and 
a diffeomorphism since $\mathbb R^n$ is simply connected. 
In each case, the set of $g$ used in the definition of $\Diff_{\mathcal A}$ is open in 
$\mathcal A(\mathbb R^n)^n$. 

Let us next check that $\on{Diff}_{\mathcal A}(\mathbb R^n)_0$ is closed under
composition.
We have  
\begin{equation}
((\Id+f)\o(\Id+g))(x)=x+g(x)+f(x+g(x)),
\tag{1}\end{equation}
and we have to check that  $x\mapsto f(x+g(x))$ is
in $\mathcal A(\mathbb R^n)$ if $f,g\in\mathcal A(\mathbb R^n)^n$.
For $\mathcal A=\mathcal B$ this follows 
by the Fa\`a di Bruno formula \ref{Faa}.
For $\mathcal A=\mathcal S$ we need furthermore the following estimate:
\begin{equation}
(\p_x^\al f)(x+g(x))=O\Bigl(\frac{1}{(1+|x+g(x)|^2)^{k}}\Bigr)
  =O\Bigl(\frac{1}{(1+|x|^2)^{k}}\Bigr)
\tag{2}\end{equation}
which holds since 
\begin{displaymath}
\frac{1+|x|^2}{1+|x+g(x)|^2}\quad\text{  is globally bounded.}
\end{displaymath}
For $\mathcal A=H^\infty$  we also need that 
\begin{align*}
\int_{\mathbb R^n} |(\p_x^\al f)(x+g(x))|^2\, dx &=  
\int_{\mathbb R^n} |(\p^\al f)(y)|^2\, \frac{dy}{|\det(\mathbb I_n + dg)((\on{Id} +g)\i (y))|} 
\\&
\le C(g) \int_{\mathbb R^n} |(\p^\al f)(y)|^2\, dy;  
\tag{3}\end{align*}
this holds since the denominator is globally bounded away from 0 since $g$ and $dg$ vanish at $\infty$ by the lemma of 
Riemann-Lebesque.
The case $\mathcal A(\mathbb R^n)=C^\infty_c(\mathbb R^n)$ is easy 
and well known. 

Let us check next that multiplication is smooth on $\Diff_{\mathcal A}(\mathbb R^n)$.
Suppose that the curves $t\mapsto \Id+f(t,\quad)$ and $t\mapsto \Id+g(t,\quad)$
are in
$C^\infty(\mathbb R,\on{Diff}_{\mathcal{A}}(\mathbb R^n))$ which means
that the functions $f,g\in C^\infty(\mathbb R^{n+1},\mathbb R^n)$ satisfy condition
$\bullet$ of either \ref{diffb}, \ref{diffh}, \ref{diffs}, or \ref{diffc}.
We have to check that $f(t,x+g(t,x))$ satisfies the same condition $\bullet$.
For this we reread the proof that composition preserves $\Diff_{\mathcal A}(\mathbb R^n)$ and pay 
attention to the further parameter $t$.

To check that the inverse $(\Id+g)\i$ is again an element in 
$\on{Diff}_{\mathcal A}(\mathbb R^n)$ for $g\in\mathcal A(\mathbb R^n)^n$, we write
$(\Id+g)\i=\Id+f$ and we have to check that $f\in \mathcal A(\mathbb R^n)^n$.  
\begin{align}
(\Id+f)\o(\Id+g)=\Id &\implies x+g(x)+f(x+g(x))=x
\notag\\&
\implies x\mapsto f(x+g(x))=-g(x)\text{  is in }\mathcal A(\mathbb R^n)^n.
\tag{4}\end{align}
We treat again first the case $\mathcal A=\mathcal B$. We know already that $\Id+g$ is a 
diffeomorphism. By definition \ref{diffb} we have 
$\det(\mathbb I_n + dg(x))\ge \ep>0$ for some $\ep$. 
This implies that 
\begin{equation*}
\|(\mathbb I_n + dg(x))\i\|_{L(\mathbb R^n,\mathbb R^n)}\quad\text{  is globally bounded,}
\tag{5}\end{equation*} 
using the inequality $\|A\i\|\le \frac{\|A\|^{n-1}}{|\det(A)|}$ for any linear 
$A:\mathbb R^n\to \mathbb R^n$. To see this, 
we write $A=US$ with $U$ orthogonal and $S=\on{diag}(s_1\ge s_2 \ge\dots\ge s_n\ge0)$.
Then $\|A\i\|=\|S\i U\i\| =\|S\i\|=1/s_n$ and 
$|\det(A)|=\det(S)=s_1.s_2\dots s_n\le s_1^{n-1}.s_n=\|A\|^{n-1}.s_n$.

Moreover, 
\begin{multline*}
(\mathbb I_n + df(x+g(x)))(\mathbb I_n + dg(x)) = \mathbb I_n 
\\
\implies
\det(\mathbb I_n + df(x+g(x)))= \det(\mathbb I_n+dg(x))\i \ge \|\mathbb I_n + dg(x)\|^{-n}\ge \et> 
0
\end{multline*}
for all $x$.
For higher derivatives
we write the Faa di Bruno formula \ref{Faa} in the following form:
\begin{align*}
&\frac{d^p(f\o (\Id+g))(x)}{p!} = 
\\&
=\on{sym}_p\Big(
\sum_{j= 1}^p \sum_{\substack{\al\in \mathbb N_{>0}^j\\ \al_1+\dots+\al_j =p}}
\frac{d^jf(x+g(x))}{j!}\Big( 
\frac{d^{\al_1}(\Id+g)(x)}{\al_1!},\dots,
\frac{d^{\al_j}(\Id+g)(x)}{\al_j!}\Big)\Big)
\\&
= \frac{d^pf(x+g(x))}{p!}\Big(\Id+dg(x),\dots,\Id+dg(x)\Big)
\\&\quad
+\on{sym}_{p}\Big(
\sum_{j= 1}^{p-1} \sum_
{\substack{\al\in \mathbb N_{>0}^j\\ \al_1+\dots+\al_j =p \\ (h_{\al_1},\dots,h_{\al_j})}}
\frac{d^jf(x+g(x))}{j!}\Big( 
\frac{d^{\al_1}h_{\al_1}(x)}{\al_1!},\dots,
\frac{d^{\al_j}h_{\al_j}(x)}{\al_j!}\Big)\Big)
\tag{6}\end{align*}
where $h_{\al_i}(x)$ is $g(x)$ for $\al_i>1$ (there is always such an $i$), and where 
$h_{\al_i}(x)=x$ or $g(x)$ if $\al_i=1$.
Now we argue as follows: 
The left hand side is globally bounded. By \thetag{5}, we know that
$\mathbb I_n + dg(x):\mathbb R^n\to \mathbb R^n$ is invertible with 
$\|(\mathbb I_n+dg(x))\i\|_{L(\mathbb R^n,\mathbb R^n)}$ globally bounded. 
Thus we can conclude by induction on $p$ that $d^pf(x+g(x))$ 
is bounded uniformly in $x$, thus also uniformly in $y=x+g(x)\in \mathbb R^n$. 
For general $\mathcal A$ we note that the left hand side is in $\mathcal A$. Since we already know 
that $f\in\mathcal B$, and since $\mathcal A$ is a $\mathcal B$-module, the last term is in 
$\mathcal A$. Thus also the first term is in $\mathcal A$, and any summand there containing at 
least one $dg(x)$ is in $\mathcal A$, so the unique summand $d^pf(x+g(x))$ is also in $\mathcal A$ as a 
function of $x$. It is thus in rapidly decreasing or in $L^2$ as a function of $y=x+g(x)$, by arguing as in \thetag{2} or 
\thetag{3} above.
Thus inversion maps $\on{Diff}_{\mathcal{A}}(\mathbb R)$ into itself. 

Next we check that inversion is smooth on $\Diff_{\mathcal A}(\mathbb R^n)$. We retrace the proof 
that inversion preserves $\Diff_{\mathcal A}$ assuming that $g(t,x)$ satisfies condition 
$\bullet$ of either \ref{diffb}, \ref{diffh}, \ref{diffs}, or \ref{diffc}.
We see again that 
$f(t,x+g(t,x))=-g(t,x)$ satisfies
the condition \ref{modelspaces} as a function of $t,x$, and we claim that $f$ then does the
same. We reread the proof paying attention to the parameter $t$ and see that the same condition 
$\bullet$ is satisfied.

We claim that $\on{Diff}_A(\mathbb R)$ is also a regular Lie group in the sense of \ref{reglie}. So let
$t\mapsto X(t,\quad)$ be a smooth curve in the Lie algebra 
$\X_{\mathcal A}(\mathbb R^n)=\mathcal{A}(\mathbb R^n)^n$, i.e., $X$ satisfies condition 
$\bullet$ of either \ref{diffb}, \ref{diffh}, \ref{diffs}, or \ref{diffc}.
The evolution of this time dependent vector field is 
the function given by the ODE
\begin{align}
&\on{Evol}(X)(t,x) = x+f(t,x),
\notag\\
&\begin{cases} \p_t (x+f(t,x)) =f_t(t,x)= X(t,x+f(t,x)),\\
f(0,x)=0. \end{cases}
\tag{7}\end{align}
We have to show first that $f(t,\quad)\in\mathcal A(\mathbb R^n)^n$ for each $t\in \mathbb R$, 
second that it is smooth in $t$ with values in $\mathcal A(\mathbb R^n)^n$, and third that 
$X\mapsto f$ is also smooth.  
For $0\le t\le C$ we consider  
\begin{equation}
|f(t,x)|\le \int_0^t |f_t(s,x)|ds=\int_0^t|X(s,x+f(s,x))|\,ds.
\tag{8}\end{equation}
Since $\mathcal A\subseteq \mathcal B$, the vector field $X(t,y)$ is uniformly bounded in 
$y\in \mathbb R^n$, locally in $t$. So the same is
true for $f(t,x)$ by \thetag{8}. 
\newline
Next consider 
\begin{align*}
\p_t d_x f(t,x) &= d_x(X(t,x+f(t,x))) 
\\&
= (d_x X)(t,x+f(t,x)) + (d_x X)(t,x+f(t,x)).d_xf(t,x)
\tag{9}\\
\|d_x f(t,x)\| &\le \int_0^t \|(d_x X)(s,x+f(s,x))\| ds 
\\&\quad
+ \int_0^t \|(d_x X)(s,x+f(s,x))\|.\|d_x f(s,x)\| ds
\\&
=: \al(t,x) + \int_0^t \be(s,x).\|d_x f(s,x)\| ds
\end{align*}
By the Bellman-Gr\"onwall inequality,  
$$
\|d_x f(t,x)\| \le 
\al(t,x) + \int_0^t \al(s,x).\be(s,x).e^{\int_s^t \be(\si,x)\,d\si}\, ds,
$$
which is globally bounded in $x$, locally in $t$.
For higher derivatives in $x$ (where $p>1$) we use Fa\'a di Bruno's formula in the form
\begin{align*}
&\p_t d_x^p f(t,x) = d_x^p(X(t,x+f(t,x))) 
=\on{sym}_p\Big(
\\&
\sum_{j= 1}^p \sum_{\substack{\al\in \mathbb N_{>0}^j\\ \al_1+\dots+\al_j =p}}
\frac{(d_x^jX)(t,x+f(t,x))}{j!}\Big( 
\frac{d_x^{\al_1}(x+f(t,x))}{\al_1!},\dots,
\frac{d_x^{\al_j}(x+f(t,x))}{\al_j!}\Big)\Big)
\\&
= (d_xX)(t,x+f(t,x))\big(d_x^p f(t,x)\big)
+\on{sym}_p\Big(
\\&
\sum_{j= 2}^{p} \sum_
{\substack{\al\in \mathbb N_{>0}^j\\ \al_1+\dots+\al_j =p }}
\frac{(d_x^jX)(t,x+f(t,x))}{j!}\Big( 
\frac{d_x^{\al_1}(x+f(t,x))}{\al_1!},\dots,
\frac{d_x^{\al_j}(x+f(t,x))}{\al_j!}\Big)\Big)
\end{align*}
We can assume recursively that $d_x^j f(t,x)$ is globally bounded in $x$, locally in $t$, for 
$j<p$. Then we have reproduced the situation of \thetag{9} (with values in the space of 
symmetric $p$-linear mappings $(\mathbb R^n)^p\to \mathbb R^n$) and we can repeat the argument 
above involving the Bellman-Gr\"onwall inequality to conclude that $d_x^pf(t,x)$ is globally bounded 
in $x$, locally in $t$. To conclude the same for $\p_t^md_x^p f(t,x)$ we just repeat  the last 
arguments for $\p_t^m f(t,x)$. So we have now proved that 
$f\in C^\infty(\mathbb R,\X_{\mathcal B}(\mathbb R^n))$. Since $x\mapsto x+f(t,x)$ is a 
diffeomorphism for each $t$ as the solution of a flow equation, it is thus in 
$\Diff_{\mathcal B}(\mathbb R^n)$.
In order to prove that 
$C^\infty(\mathbb R,\X_{\mathcal B}(\mathbb R^n))\ni X\mapsto \on{Evol}(X)(1,\quad)\in \Diff_{\mathcal B}(\mathbb R^n)$
is smooth, we consider a smooth curve $X$ in $C^\infty(\mathbb R,\X_{\mathcal B}(\mathbb R^n))$; 
thus $X(t_1,t_2,x)$ is smooth on $\mathbb R^2\x \mathbb R^n$, globally bounded in $x$ in each 
derivative separately, locally in $t=(t_1,t_2)$ in each derivative. Or, we assume that $t$ is 
2-dimensional in the argument above. But then it suffices to show that 
$(t_1,t_2)\mapsto X(t_1,t_2,\quad)\in \X_{\mathcal B}(\mathbb R^n)$ is smooth along smooth curves in 
$\mathbb R^2$, and we are again in the situation we have just treated. 
Thus $\Diff_{\mathcal B}(\mathbb R^n)$ is a regular Lie group. 

If $\mathcal A=\mathcal S$, we already know that $f(s,x)$ is globally bounded in $x$, locally in 
$s$. Thus we may insert 
$X(s,x+f(s,x))=O(\frac1{(1+|x+f(s,x)|^2)^k})=O(\frac1{(1+|x|^2)^k})$ into 
\thetag{8} and can conclude that $f(t,x)=O(\frac1{(1+|x|^2)^k})$
globally in $x$, locally in $t$, for each $k$.
Using this argument, we can repeat the proof for the case $\mathcal A=\mathcal B$ from above and 
conclude that $\Diff_{\mathcal S}(\mathbb R^n)$ is a regular Lie group. 
\newline
If $\mathcal A= H^\infty$ we first consider the differential version of \thetag{8}, 
\begin{align*}
\|d_x f(s,x)\| &= 
 \Big\|\int_0^t d_x(X(s,\quad))(x+f(s,x)).(\mathbb I_n + df(s,\quad)(x))\, ds \Big\|
\\&
\le  \int_0^t \big\|d_x(X(s,\quad))(x+f(s,x))\big\|. C\, ds 
\tag{10}\end{align*}
since $d_x f(s,x)$ is globally 
bounded in $x$, locally in $s$, by the case $\mathcal A =\mathcal B$. The same holds for $f(s,x)$. 
Moreover, $X(s,x)$ vanishes at $x=\infty$ by the lemma of 
Riemann-Lebesque for each $x$ and it is continuous in all variables, so that the same holds for 
$f(s,x)$ by \thetag{8}.
Now we consider
\begin{align*}
\int_{\mathbb R^n}\|(d_x^p f)(t,x)\|^2\,dx 
&= \int_{\mathbb R^n}\Big\|\int_0^t d_x^p\big(X(s,\Id+f(s,\quad))\big)(x)\,ds\Big\|^2\,dx.
\tag{11}\end{align*}
We apply the Fa\'a di Bruno formula in the form \thetag{6} to the integrand,
remember that we already know that each $d^{\al_i}(\Id+ f(s,\quad))(x)$ is globally bounded, 
locally in $s$, thus \thetag{11} is		 
\begin{align*}
&\le \int_{\mathbb R^n}\Big(\int_0^t \sum_{j=1}^p\|(d_x^j X)(s,x+f(s,x))\|.C_j\,ds\Big)^2\,dx
\\&
= \int_{\mathbb R^n}\Big(\int_0^t \sum_{j=1}^p\|(d_x^j 
X)(s,y)\|.C_j\,ds\Big)^2\,\frac{dy}{|\det(\mathbb I_n + df(s,\quad))((\mathbb I_n+ f(s,\quad))\i(y)| }
\end{align*}
which is finite since $X(s,\quad)\in H^\infty$ and since the determinand in the denominator is 
bounded away from zero -- we just checked that $d_xf(s,\quad)$ vanishes at infinity.
Then we repeat this for $\p_t^m d_x^p f(t,x)$.
This shows that $\on{Evol}(X)(t,\quad)\in \Id + H^\infty(\mathbb R^n)^n$. As solution of an 
evolution equation for a bounded non-autonomous vector field it is a diffeomorphism, and thus in 
$ \Diff_{H^\infty}(\mathbb R^n)$ for each $t$. By the 
same trick as in the case $\mathcal A=\mathcal B$ we can conclude that 
$\Diff_{H^\infty}(\mathbb R^n)$ is a regular Lie group. 

We prove now that $\Diff_{\mathcal S}(\mathbb R^n)$ is a normal subgroup of 
$\Diff_{\mathcal B}(\mathbb R^n)$. So let $g\in \mathcal B(\mathbb R^n)^n$ with 
$\det(\mathbb I_n + dg(x))\ge \ep>0$ 
for all $x$, and 
$s\in \mathcal S(\mathbb R^n)^n$ with $\det(\mathbb I_n + ds(x)) >0$ for all $x$.
We consider  
\begin{align*}
(\Id +g)\i(x) &= x + f(x)\quad\text{  for } f\in \mathcal B(\mathbb R^n)^n \iff f(x+g(x))=-g(x)
\\
((\Id+g)\i\o(&\Id+s)\o (\Id+g)\big)(x) = ((\Id+f)\o(\Id+s)\o(\Id+g)\big)(x) = 
\\&
= x + g(x) + s(x + g(x)) + f\big(x + g(x) + s(x + g(x))\big)
\\&
= x + s(x + g(x)) - f(x+g(x)) + f\big(x + g(x) + s(x + g(x))\big).
\end{align*}
Since $g(x)$ is globally bounded we get 
$s(x+g(x)) = O((1 + |x + g(x)|)^{-k}) = O((1 + |x|)^{-k})$ for each $k$. 
For $d_x^p(s\o (\Id + g))(x)$ this follows from the Fa\'a die Bruno formula in the form of 
\thetag{6}.
Moreover we have
\begin{multline*}
 f\big(x + g(x) + s(x + g(x))\big) - f(x+g(x)) =
 \\
= \int_0^1 df\big(x+g(x)+ ts(x+g(x))\big)(s(x+g(x)))\,dt
\end{multline*}
which is in $\mathcal S(\mathbb R^n)^n$ as a function of $x$ since $df$ is in $\mathcal B$ and 
$s(x+g(x))$ is in $\mathcal S$. 

Finally we prove that $\Diff_{H^\infty}(\mathbb R^n)$ is a normal subgroup of 
$\Diff_{\mathcal B}(\mathbb R^n)$. We redo the last proof under the assumption that 
$s\in H^\infty(\mathbb R^n)^n$. By the argument in \thetag{3} we see that 
$s(x+g(x))$ is in $H^\infty$ as a function of $x$. The rest is as above since $H^\infty$ is an 
ideal in $\mathcal B$ as noted in \ref{ideal}.
\end{proof}

\begin{cor}
$\Diff_{\mathcal B}(\mathbb R^n)$ acts on $\Ga_c$, $\Ga_{\mathcal S}$ and $\Ga_{H^\infty}$ of any 
tensorbundle over $\mathbb R^n$ by pullback. The infinitesimal action of the Lie algebra $\X_{\mathcal B}(\mathbb R^n)$ on 
these spaces by the Lie derivative thus maps each of these spaces into itself. 
A fortiori, $\Diff_{H^\infty}(\mathbb R^n)$ acts on $\Ga_{\mathcal S}$ of any tensor bundle by pullback. 
\end{cor}

\demo{Proof}
Since $\Diff_c(\mathbb R^n)$, $\Diff_{\mathcal S}(\mathbb R^n)$, and 
$\Diff_{H^\infty}(\mathbb R^n)$ are normal subgroups in $\Diff_{\mathcal B}(\mathbb R^n)$, their Lie algebras  
$\X_{\mathcal A}(\mathbb R^n)=\Ga_{\mathcal A}(T\mathbb R^n)$ are all invariant under the adjoint 
action of $\Diff_{\mathcal B}(\mathbb R^n)$. This extends to all tensor bundles. The Lie 
derivatives are just the infinitesimal versions of the adjoint actions. 
\qed\enddemo

\bibliographystyle{plain}
\def\cprime{$'$}

\end{document}